\documentclass{lms}
\usepackage[latin1]{inputenc}
\usepackage{listings}
\usepackage{verbatim}
\usepackage{color}
\usepackage{textcomp, latexsym}
\usepackage{graphicx}
\usepackage{float}
\usepackage{flafter}
\usepackage{longtable}
\usepackage{amsmath}
\usepackage{amsfonts}
\usepackage{amssymb}
\usepackage[bookmarks=false]{hyperref}
\usepackage{epsfig}
\usepackage{color}
\usepackage{pstricks,pst-node,pst-text,pst-3d}

\newnumbered{definition}{Definition}
\newnumbered{remark}{Remark}
\newnumbered{algorithm}{Algorithm}
\newnumbered{example}{Example}
\title[Discrete Variational Calculus]{Discrete Variational Calculus for B-spline Approximated Curves}
\author{Jun Zhao, Elizabeth Mansfield}
\date{\today}
\classno{11B83, 49M25, 42A10, 42A15, 70G75}
\begin{document}
\maketitle
\vspace{-0.5cm}\begin{center}{\em School of Mathematics, Statistics \& Actuarial Science}\end{center}
\vspace{0cm}\begin{center}{\em Cornwallis Building}\end{center}
\vspace{0.03cm}\begin{center}{\em University of Kent}\end{center}
\vspace{0.06cm}\begin{center}{\em Canterbury}\end{center}
\vspace{0.09cm}\begin{center}{\em Kent CT2 7NF}\end{center}
\vspace{0.12cm}\begin{center}{jz45@kent.ac.uk,\quad e.l.mansfield@kent.ac.uk}\end{center}
\vspace{0.5cm}
\begin{abstract}
We study variational problems for curves approximated by B-spline curves. We show that one can obtain discrete Euler-Lagrange equations for the data describing the approximated curves.

Our main application is to the occluded curve problem in 2D and 3D. In this case, the aim is to find various aesthetically pleasing solutions as opposed to a solution of a physical problem. The Lagrangians of interest are invariant under the special Euclidean group action for which B-spline approximated curves are well suited.

Smooth Lagrangians with special Euclidean symmetries involve curvature, torsion, and arc length. Expressions in these, in the original coordinates, are highly complex. We show that, by contrast, relatively simple discrete Lagrangians offer excellent results for the occluded curve problem. The methods we develop for the discrete occluded curve problem are general and can be used to solve other discrete variational problems for B-spline curves. 
\end{abstract}
\section{Introduction}
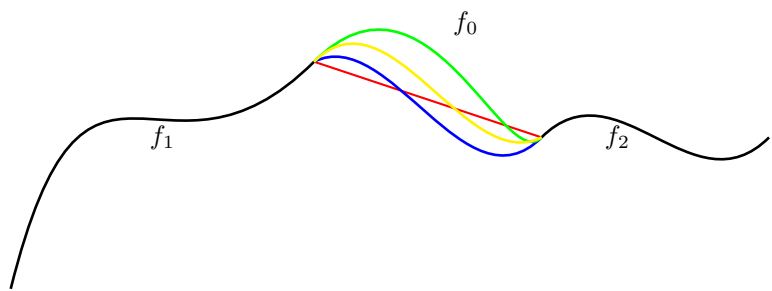
\begin{figure}[H]
\begin{pspicture}(-3,4)
\psbezier[linewidth=1pt,showpoints=false]{}(0,0)(1,4)(2,1)(4,3)
\psbezier[linewidth=1pt,showpoints=false]{}(7,2)(8,3)(9,1)(10,2)
\psline[linecolor=red](4,3)(7,2)
\psbezier[linewidth=1pt,showpoints=false,linecolor=blue]{}(4,3)(5,3.5)(6,1)(7,2)
\psbezier[linewidth=1pt,showpoints=false,linecolor=green]{}(4,3)(5.5,4.5)(6.5,1.5)(7,2)
\psbezier[linewidth=1pt,showpoints=false,linecolor=yellow]{}(4,3)(5,4)(6,1.5)(7,2)
\rput(2,2){$f_1$}
\rput(8,2){$f_2$}
\rput(6,3.5){$f_0$}
\end{pspicture}
  \caption{Determine a path $f_0$ connecting the paths $f_1$ and $f_2$}\label{path}
\end{figure}

Finding a path which connects two existing paths as described in Figure \ref{path} is a practical problem which has applications in physics, engineering, computer vision and many other areas; mathematically, the occluded curve problem leads to the study of the Euler-Lagrange equations derived from the extremization of a variational problem.

As we can see from Figure \ref{path}, there are many possible solutions, even a straight line could be a ``solution'' of the occluded curve. Therefore, certain constrains have to be applied to the solution curve:

\begin{itemize}
\item The solution curve has to be smooth, which means the boundary conditions must ensure that the occluded curve has $C^1$ continuity at the boundary;

\item The solution method needs to be equivariant under certain transformation (such as rotation and translation), which means that completing the curve before translation and rotation is the same as completing the curve after;

\item The solution curve needs to satisfy certain optimization criteria.

\end{itemize}
For plane curves, these three conditions can be satisfied simultaneously if the solution curve $f(s)=(x(s),y(s))$ extremizes a functional $\int L(s,\kappa,\kappa_s\cdots){\rm d}s$, where $L{\rm d}s$ is the Lagrangian, and $L$ is a function of the invariants of the Euclidean group, where $s$ is the arc length, and $\kappa$ is the Euclidean  curvature of $f(s)$. For space curves, torsion and its derivatives may be included in the Lagrangian. 

The solution curve then satisfies the Euler-Lagrange equations for this variational problem. The method of solving the occluded curve in the smooth case can be found in the book\cite{mansfield}. We give only a brief review of the method in the next paragraph.
 
In the smooth case, the Euler-Lagrange equations are very complex. In 2D, for $L=\kappa^2$, the solution of the Euler-Lagrange equations which is known as Euler's elastica yields that the curvature $\kappa$ is a Jacobean Sine function; in 3D, for the same Lagrangian as in 2D, $L=\kappa^2$, the solutions of the Euler-Lagrange equations show that the curvature and the torsion are Weierstrass $\wp$ functions. Finding the explicit solutions such as $x(s)$ and $y(s)$ involves yet more complicated calculation. So using approximation method becomes important and necessary.

In computer vision and many other fields, all geometric entities are represented by a variety of approximation methods. We use B-spline curves in this paper to approximate given curves and develop a discrete variational calculus for the occluded curve problem based on them.

The input for the occluded curve problem is the outside pieces of the curve $f_1$ and $f_2$ (see Figure \ref{path}). We can use the control points of their B-spline approximations as discrete input data for the boundary conditions. Then we can construct the discrete Euler-Lagrange equations, which are for the control points of the missing piece of the curve.

Similarly to the smooth case, the solutions of the discrete Euler-Lagrange equations depend highly on the choice of the Lagrangian. Choosing an appropriate Lagrangian, which may or may not represents a corresponding smooth Lagrangian, we will be able to obtain an aesthetically pleasing result.

Different approaches of the similar problems are studied by others. Some used discrete points on the curve to approximate and construct invariants\cite{shakiban} which has disadvantage when dealing with complicated curves, the data need to be processed are massive. Also precision is lost when reducing points to simplify the problem.

The Euler-Lagrange equations for many kinds of discrete variational problems are difference equations. Our Euler-Lagrange equations are different. Difference equations are special type of recurrence equations which recursively define sequences; our discrete equations, on the other hand, define quantities without recurrence. The analysis and calculation of our discrete equations are much easier than difference equations and hence grant us advantages in solving the discrete occluded curve problem.

In this paper firstly we give the context of the discrete occluded curve problem. Then a brief introduction to B-spline curves and our main tool, syzygies on the control points, which we will define. Using discrete syzygies is essential to constructing the Euler-Lagrange equations. After that, we introduce the boundary conditions which simplify the discrete system by eliminating some unknowns. The simplified system can be then solved. Some examples both in 2D and 3D are given, and readers can compare the results from using different Lagrangians and input data. From those examples we can see that choosing an appropriate Lagrangian is very important of obtaining an aesthetically pleasing solution. A general algorithm is then described.

By considering different sorts of  boundary conditions, more features can be introduced into the solution curve. Examples are given in 2D space which are classified into different cases based on the boundary conditions and input data. Finally, a revised algorithm is given which is designed not only for the discrete occluded curve problem but also for the discrete variational problem in general.

Without further ado, let us now introduce the discrete occluded curve problem.
\section{Discrete Occluded Curve Problem}
The discrete occluded curve problem is based on approximation. In the industrial field of computer graphic design, using polynomials to approximate actual curves is very practical, and B-spline\cite{deboor} has proved to be one of the most efficient methods due to its outstanding properties. We will use B-spline curves to simplify our variational problem especially because of its affine invariance property, which implies, that the curve after transformation can be re-constructed by the transformed control points. Hence we can convert the differential invariants problem to the finite difference invariants problem, which is much easier.

Talking about invariants, we need to give the group actions which act on our curves. And the geometric entities such as curvature and torsion which will be used to analyze the variational problem are invariant under those group actions. The group actions we will use in this paper are those of the special Euclidean group actions ${\rm SE}(n)$, which include rotation and translation in $n$ dimensional space.

\textbf{Notation:} {\em I will use $\bullet$ to represent the group action, i.e. $\widetilde{x}=g\bullet x,\;g\in{\rm SE}(n)$; and $\cdot$ to represent the inner product between two vectors.}

Considering the input curves are at the ``normalized'' positions, instead of finding the solution curve, we only need to find the control points of the solution curve. This leads us to the problem of finding and solving discrete Euler-Lagrange equations.
\subsection{B-spline Curves and Syzygies}
Before we start the discussion of the discrete Euler-Lagrange equations, we will give a brief overview of the B-spline basis functions and syzygies\cite{mansfield}.

In computer graphics and computer aided design area, the B-spline usually refers to the parametric spline curves that are constructed by interpolating the control points using the B-spline basis functions which are defined over a certain knot sequence. The control points of a B-spline curve, just as the name implied, are the points which control the shape of the spline curve. The basis functions of a B-spline are defined recursively over a knot sequence, which is a sequence of monotone increasing real numbers. When the knots are equidistant, the spline is called {\it uniform}, otherwise {\it non-uniform}. The degree of the spline curve depends on the degree of the basis functions.

Let a B-spline curve be given as $$f_j(t)=B(t)P=\sum_{i=j}^nb^j_{i,k}(t)p_i,\;t\in[t_{j+k},t_{j+k+1}],\quad j=1..n-k$$ where $B(t)=(b^j_{i,k}(t))$ is a row vector of basis functions which are defined over the knot sequence $\{t_1,t_2,\cdots,t_m\}$, $P=(p_i)$ is a column vector of control points, $j$ is sequential number of the spline curve, and $k$ is the degree of the spline.

\textbf{Notation:} {\em When $n=k+1$, it is a one-piece spline curve, where $j=1$, the subscript of $f_j(t)$ will be omitted; when $n>k+1$, it is a multi-piece spline curve, where $j=1..n-k$. And $m..n$ as in Maple means a set of integers from $m$ to $n$}

The basis functions are defined recursively over a knot sequence $\{t_1,t_2,\cdots,t_m\}$, with $0\leq t_i\leq t_{i+1}\leq1,\forall i$, $m=n+k+1$, where $n$ is the number of the control points:
\begin{equation}\begin{array}{lll}\label{basis}
b^j_{i,0}(t)&=&\left\{\begin{array}{ll}1& t_i\leq t<t_{i+1}\;\&\;t_i<t_{i+1}\\ 0& otherwise\end{array}\right.\\[20pt]
b^j_{i,k}(t)&=&\displaystyle\frac{t-t_i}{t_{i+k}-t_i}b^j_{i,k-1}(t)+\displaystyle\frac{t_{i+k+1}-t}{t_{i+k+1}-t_{i+1}}b^j_{i+1,k-1}(t)
\end{array}\end{equation}
Let a B-spline curve $f(t),\;t\in[0,1]$ with control points $(p_i)_{i=1}^n$ be given. In order to have certain properties of the B-spline curve, such as $f(0)=p_1,f(1)=p_n$, line segments $\overrightarrow{p_1p_2}$ and $\overrightarrow{p_{n-1}p_n}$ are tangent vector at $f(0)$ and $f(1)$ respectively (\emph{Bessel end conditions}), we need to set the first $k+1$ knots to 0 and the last $k+1$ knots to 1. In the following analysis and examples we use the uniform B-spline curves which satisfy the Bessel end conditions.

Considering the control points of the B-spline curves at the normalized positions as discrete invariants, the differences between the control points are difference invariants. Those difference invariants will be utilized in the Lagrangian.

In the variational problem using discrete Lagrangians, in order to find the extremized solution, we introduce a dummy variable $\epsilon$, and apply a variation to the difference invariants with respect to $\epsilon$. The syzygies are the relations between the variation of different difference invariants.

More definition of syzygies in smooth case can be found in the book\cite{mansfield}. Following sections will give the construction of the syzygies in the discrete case.
\subsection{Construction of Discrete Euler-Lagrange Equations}
Considering now the control points of the given curves as input data. The discrete Lagrangian is defined over the whole interval of control points, but the variation only applies to the control points of the occluded curve. The output will be the control points of occluded curves, which extremize the variational problem.

Firstly we discuss the construction of the discrete Euler-Lagrange equations bases on single-piece and then multi-piece cubic B-spline curves which satisfy the Bessel end conditions respectively.
\subsubsection{One-Piece B-spline plane Curve}
Let a one-piece cubic B-spline plane curve be $f(t)=\sum_{i=1}^4b_{i,3}(t)p_i$, where the basis functions are defined over the knot sequence $\{0,0,0,0,1,1,1,1\}$, be a solution to the Euler-Lagrange equations derived from the extremization of the variational problem (\ref{variation1}). The corresponding Lagrangian is invariant under the group action $g\in{\rm SE}(2)$. We have
\begin{equation}\label{variation1}
\mathcal L[I,I_K]=\displaystyle\sum_iL(I_i,I_{i,K}),
\end{equation}
where the invariants $I=(I_i)=(p_i)$ are the control points at the normalized position, $I_K=(I_{i,K})=(p^K_i)$ are the differences of control points at the normalized position. Here $K=\{1,2,3\}$ is a set of difference orders. The normalized position is defined as the last control point of left input curve locating at the origin, and the first control point of right input curve locating at the $x$-axis (refer to Figure \ref{path}). The difference of the control points are defined as:
$$p_i^l=p_{i+1}^{l-1}-p_i^{l-1},\quad p_i^0=p_i,\quad l=1..3$$
Now introduce a dummy variable $\epsilon$ and apply a variation to the difference invariants with the respect to $\epsilon$, so that $I_{i}^\epsilon,\; i=1..4$ is the variation of $I_i$, and the syzygies of the difference invariants for the control points are:
\begin{equation}\label{discretesyzygy}\begin{array}{rcl}
I_{i,1}&=&I_{i+1}-I_i=(S-{\rm id})I_i\\[8pt]
I_{i,2}&=&I_{i+1,1}-I_{i,1}=(S-{\rm id})I_{i,1}=(S-{\rm id})^2I_i\\[8pt]
&\cdots&\\[8pt]
\displaystyle\frac{{\rm d}}{{\rm d}\epsilon}I_i&=&I_{i}^\epsilon\\[8pt]
\displaystyle\frac{{\rm d}}{{\rm d}\epsilon}I_{i,1}&=&I_{i+1}^\epsilon-I_{i}^\epsilon=(S-{\rm id})I_{i}^\epsilon\\[8pt]
\displaystyle\frac{{\rm d}}{{\rm d}\epsilon}I_{i,2}&=&\displaystyle\frac{{\rm d}}{{\rm d}\epsilon}(I_{i+1,1}-I_{i,1})=(S-{\rm id})^2I_{i}^\epsilon\\[8pt]
&\cdots&
\end{array}\end{equation}
where $S$ is the shift map, $S(I_i)=I_{i+1}$, and ${\rm id}$ is the identity map. Notice that the shift operator and the group action do not commute, which means $S(I_{i,1})=S(g\bullet (p_{i+1}-p_i))\neq g\bullet(S(p_{i+1}-p_i))$. The shift operator here only is defined on the normalized difference of the control points. Since the variation only happens on the occluded curve, the shift operator and its inverse operator will send the variation $I_{i}^\epsilon,\; i\geq n\,\&\,i\leq 1$ to zero, corresponding to the input data being fixed. Hence the extremization of the variational problem (\ref{variation1}) becomes:
\begin{equation}\label{discrete}\begin{array}{rcl}
\left.\displaystyle\frac{{\rm d}}{{\rm d}\epsilon}\right|_{\epsilon=0}\mathcal L&=&\left.\displaystyle\frac{{\rm d}}{{\rm d}\epsilon}\right|_{\epsilon=0}\displaystyle\sum_iL(I_i,I_{i,K})\\[8pt]
0&=&\displaystyle\sum_i\left(\displaystyle\frac{{\rm d}L}{{\rm d}I_i}(\displaystyle\frac{{\rm d}}{{\rm d}\epsilon}I_i)+\displaystyle\sum_{j=1}^{3}\displaystyle\frac{{\rm d}L}{{\rm d}I_{i,K_j}}(\displaystyle\frac{{\rm d}}{{\rm d}\epsilon}I_{i,K_j})\right)\\[8pt]
&=&\displaystyle\sum_i\left(\displaystyle\frac{{\rm d}L}{{\rm d}I_i}I_{i}^\epsilon+\displaystyle\sum_{j=1}^{3}\displaystyle\frac{{\rm d}L}{{\rm d}I_{i,K_j}}(S-{\rm id})^{j}I_{i}^\epsilon\right)\\[8pt]
&=&\displaystyle\sum_i\left(\left(\displaystyle\frac{{\rm d}L}{{\rm d}I_i}+\displaystyle\sum_{j=1}^{3}(S^{-1}-{\rm id})^{j}\displaystyle\frac{{\rm d}L}{{\rm d}I_{i,K_j}}\right)I_{i}^\epsilon+{\rm B.T.s}\right)\\[8pt]
&=&\displaystyle\sum_i\left(E_i(L)I_{i}^\epsilon+{\rm Boundary\;Terms}\right)
\end{array}\end{equation}
The discrete Lagrangian does not need to be exact discretization of smooth Lagrangian. It just needs to be invariant and have some geometric meanings, easily solved, and the solution should be aesthetically pleasing. Now let us see some examples.
\begin{example}\label{example1}
Let an occluded curve be one-piece cubic B-spline curve, with the basis functions defined over the knot sequence $\{0,0,0,0,1,1,1,1\}$ and the control points $p_i, i=1..4$. Choose the Lagrangian $L=I_{1,2}\cdot I_{2,2}$, which represents the angle between two second order differences of the control points, with $I_{1,2}=(S-{\rm id})^2I_1$ and $I_{2,2}=(S-{\rm id})^2I_2$. Then the corresponding Euler-Lagrange equations are:
\begin{equation}\label{discrete1}E_1=(S^{-1}-{\rm id})^2\left[(1,1)^T\cdot I_{2,2}\right]\quad E_2=(S^{-1}-{\rm id})^2\left[I_{1,2}\cdot(1,1)^T\right]\end{equation}
\end{example}
\begin{example}\label{example2}
Let an occluded curve be the same as in Example \ref{example1}, now let the Lagrangian be $L=I_{1,1}\cdot I_{3,1}$, which represents the angle between two first order differences of the control points, with $I_{1,1}=(S-{\rm id})I_1$ and $I_{3,1}=(S-{\rm id})I_3$. Then the corresponding Euler-Lagrange equations are:
\begin{equation}\label{discrete2}E_1=(S^{-1}-{\rm id})\left[(1,1)^T\cdot I_{3,1}\right]\quad E_2=(S^{-1}-{\rm id})\left[I_{1,1}\cdot(1,1)^T\right]\end{equation}
\end{example}
As we can see from the examples above, choosing different Lagrangians will give us different Euler-Lagrange equations, thus different solutions to the occluded curve.

The basis functions of a cubic B-Spline curve are already given, so there are 4 control points which contain eight unknowns to determine, but we only got two Euler-Lagrange equations. Hence we need some boundary conditions to eliminate the extra unknowns.\newline

\textbf{Determine the boundary conditions}

\begin{figure}[H]
\begin{pspicture}(-2.3,4)
\psbezier[linewidth=1pt,showpoints=true]{}(0,0)(1,4)(2,1)(4,3)
\psbezier[linewidth=1pt,showpoints=true]{}(7,2)(8,3)(9,1)(10,2)
\psbezier[linewidth=1pt,showpoints=true,linecolor=red]{}(4,3)(5,4)(6,1)(7,2)
\rput(2,2){$f_1$}
\rput(8,2){$f_2$}
\rput(6,3.5){$f_0$}
\rput(3.8,3.3){$p_1$}
\rput(4.8,4.3){$p_2$}
\rput(6.1,0.8){$p_3$}
\rput(6.8,2.3){$p_4$}
\rput(8,3.3){$p_5$}
\rput(9,0.8){$p_6$}
\rput(10.1,2.3){$p_7$}
\rput(2,0.8){$p_{0}$}
\rput(1,4.3){$p_{-1}$}
\rput(0,-0.2){$p_{-2}$}
\end{pspicture}
  \caption{Determine the boundary conditions of the solution path $f_0$}\label{path1}
\end{figure}
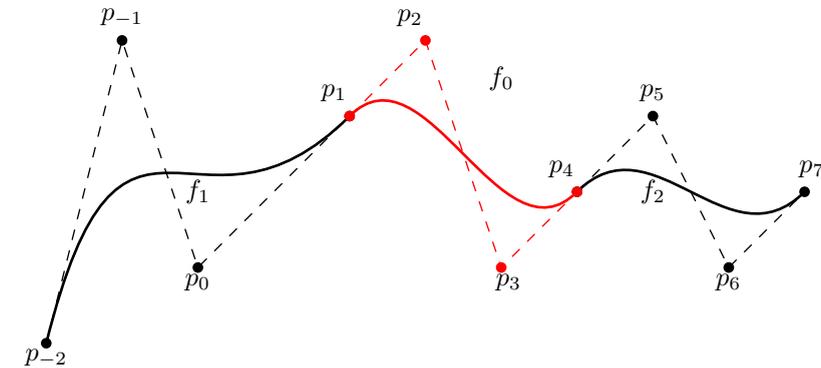
Figure \ref{path1} shows the occluded curve problem with the control points $p_i=(x_i,y_i)^T$. As we can see the control points (red) of the solution curve must satisfy some boundary conditions:

1. The first and the last control points ($p_1$ and $p_4$) of the occluded curve are the last and the first control points of the two input curves ($f_1$ and $f_2$) respectively;

2. The second and the second last control points ($p_2$ and $p_3$) must lie on the tangent line of the first and last curve points respectively.

According to the two rules above, if we set 
\begin{equation}\label{boundary}p_2=\alpha(p_1-p_0)\qquad p_3=\beta(p_4-p_5),\end{equation}
with $\alpha,\beta\in\mathbb R$, as shown in Figure \ref{path1}, we can reduce the eight unknowns to two, thus we can use the two Euler-Lagrange equations in (\ref{discrete1}) and (\ref{discrete2}) to solve the control points of the occluded curve.

Since we require the solution curve are calculated at the normalized positon, we can transform the original input curves to the normalized position first and calculate the solution curve, then transform them all back to the original position.\newline

\textbf{Solving the discrete Euler-Lagrange equations}\newline

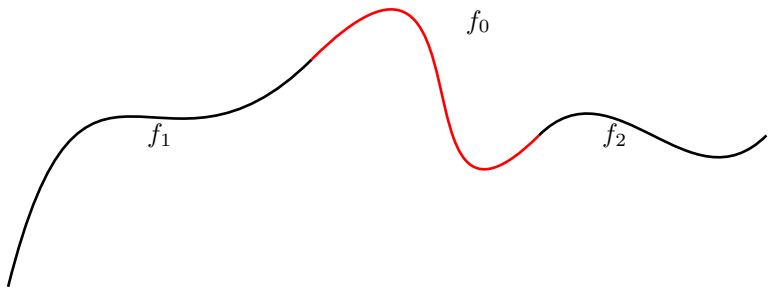
\begin{figure}[H]
\begin{pspicture}(4,4)
\psbezier[linewidth=1pt,showpoints=false]{}(0,0)(1,4)(2,1)(4,3)
\psbezier[linewidth=1pt,showpoints=false]{}(7,2)(8,3)(9,1)(10,2)
\psbezier[linewidth=1pt,showpoints=false,linecolor=red]{}(4,3)(6.5,5.5)(5,0)(7,2)
\rput(2,2){$f_1$}
\rput(8,2){$f_2$}
\rput(6.2,3.5){$f_0$}
\end{pspicture}
  \caption{Solution curve $f_0$ to the problem in Figure \ref{path1} using Euler-Lagrange equations in (\ref{discrete1}) is shown in red}\label{solution}
\end{figure}

{\em \textbf{Solving Example \ref{example1}}}
{\em Now let us solve the problem described in Figure \ref{path1}. Let two one-piece cubic B-spline curves $f_1$ and $f_2$ been given, their control points are $p_{-2}=(0,0), p_{-1}=(1,4),p_{0}=(2,1),p_{1}=(4,3),p_4=(7,2),p_5=(8,3),p_6=(9,1),p_7=(10,2)$. First transform them to the normalized position, that is $\widetilde{p_1}=(0,0),\widetilde{y_4}=0$, after solve the problem, transform them back to the original position. Using the Euler-Lagrange equations in (\ref{discrete1}) and the boundary conditions discussed previously, we will obtain the solution curve as in Figure \ref{solution}}:\newline

\hspace*{-5.5mm}{\em \textbf{Solving Example \ref{example2}}}
{\em Using the same input data and conditions as in Example \ref{example1}, but the Euler-Lagrange equations in (\ref{discrete2}), the corresponding solution curve is shown in Figure \ref{solution2}
\begin{figure}[H]
\begin{pspicture}(4,4)
\psbezier[linewidth=1pt,showpoints=false]{}(0,0)(1,4)(2,1)(4,3)
\psbezier[linewidth=1pt,showpoints=false]{}(7,2)(8,3)(9,1)(10,2)
\psbezier[linewidth=1pt,showpoints=false,linecolor=red]{}(4,3)(6,5)(6.8,1.8)(7,2)
\rput(2,2){$f_1$}
\rput(8,2){$f_2$}
\rput(6.2,3.5){$f_0$}
\end{pspicture}
  \caption{Solution curve $f_0$ to the problem in Figure \ref{path1} using Euler-Lagrange equations in (\ref{discrete2}) is shown in red. The solution curve is smooth, but has very high curvature at the right hand end.}\label{solution2}
\end{figure}
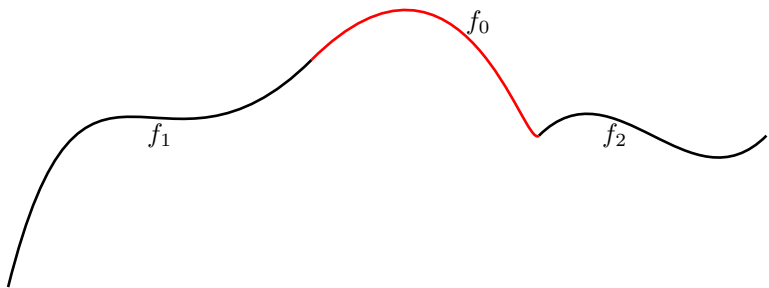
}
\begin{remark}
Comparing the two solution curves using different Lagrangians, we can see that, the solution in Figure \ref{solution} obviously is better than in Figure \ref{solution2}. That is because the Lagrangian we used in Figure \ref{solution} includes the second order difference of the control points which relates to the dominant derivatives in the curvature form. However, the Lagrangian we used in Figure \ref{solution2} only contains the first order difference, which is aesthetically less pleasing. Hence the choice of the Lagrangian is crucial in solving the occluded curve problem.
\end{remark}
Next We consider the same problem but in three dimensional space.
\subsubsection{Multi-piece Space Curve}
The construction of the basis functions of a space curve is the same as of a plane curve, we only need to add one more dimension to the control points.

\begin{figure}[H]
  \centerline{\includegraphics[width=5cm]{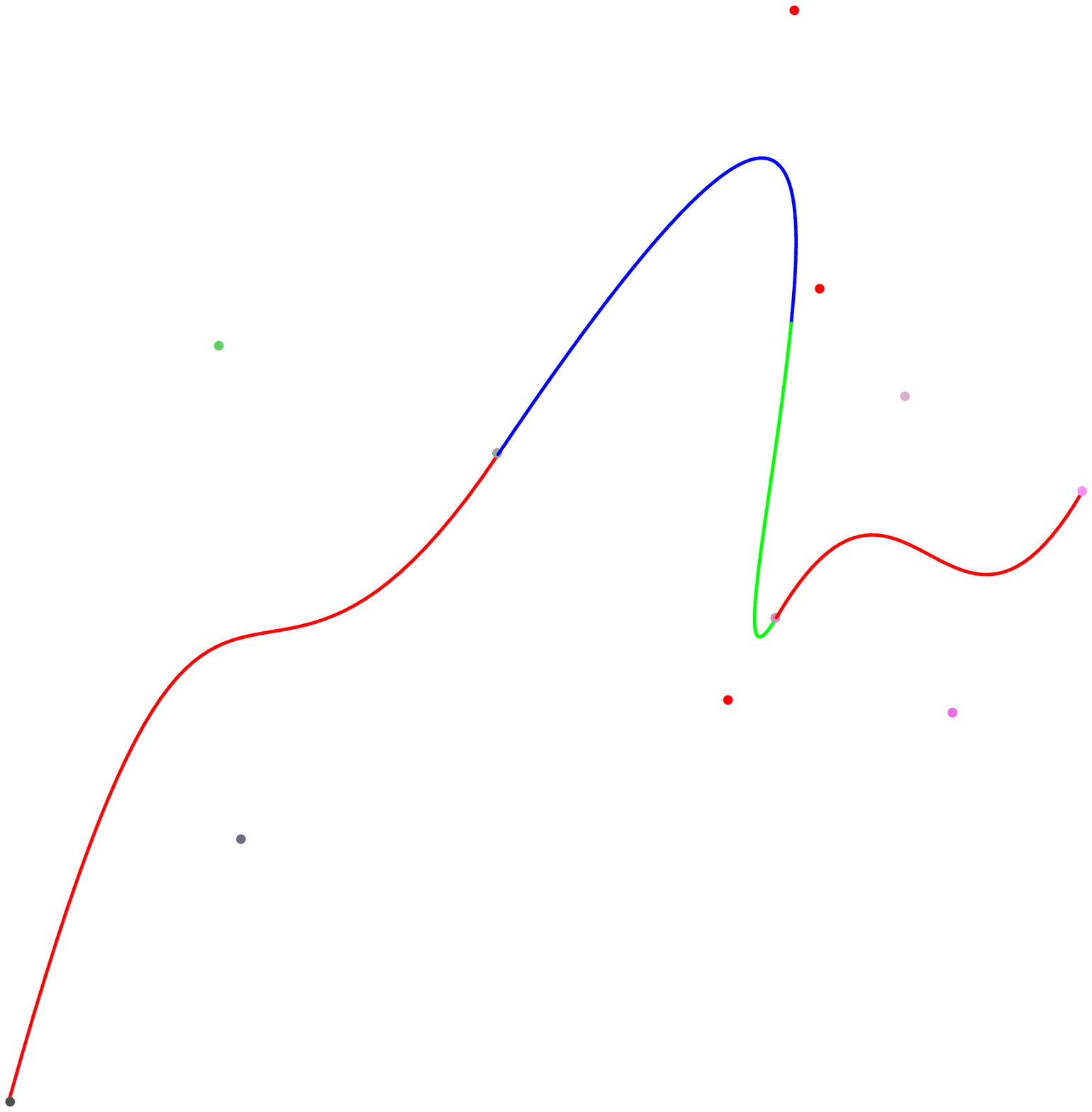}\qquad\includegraphics[width=5cm]{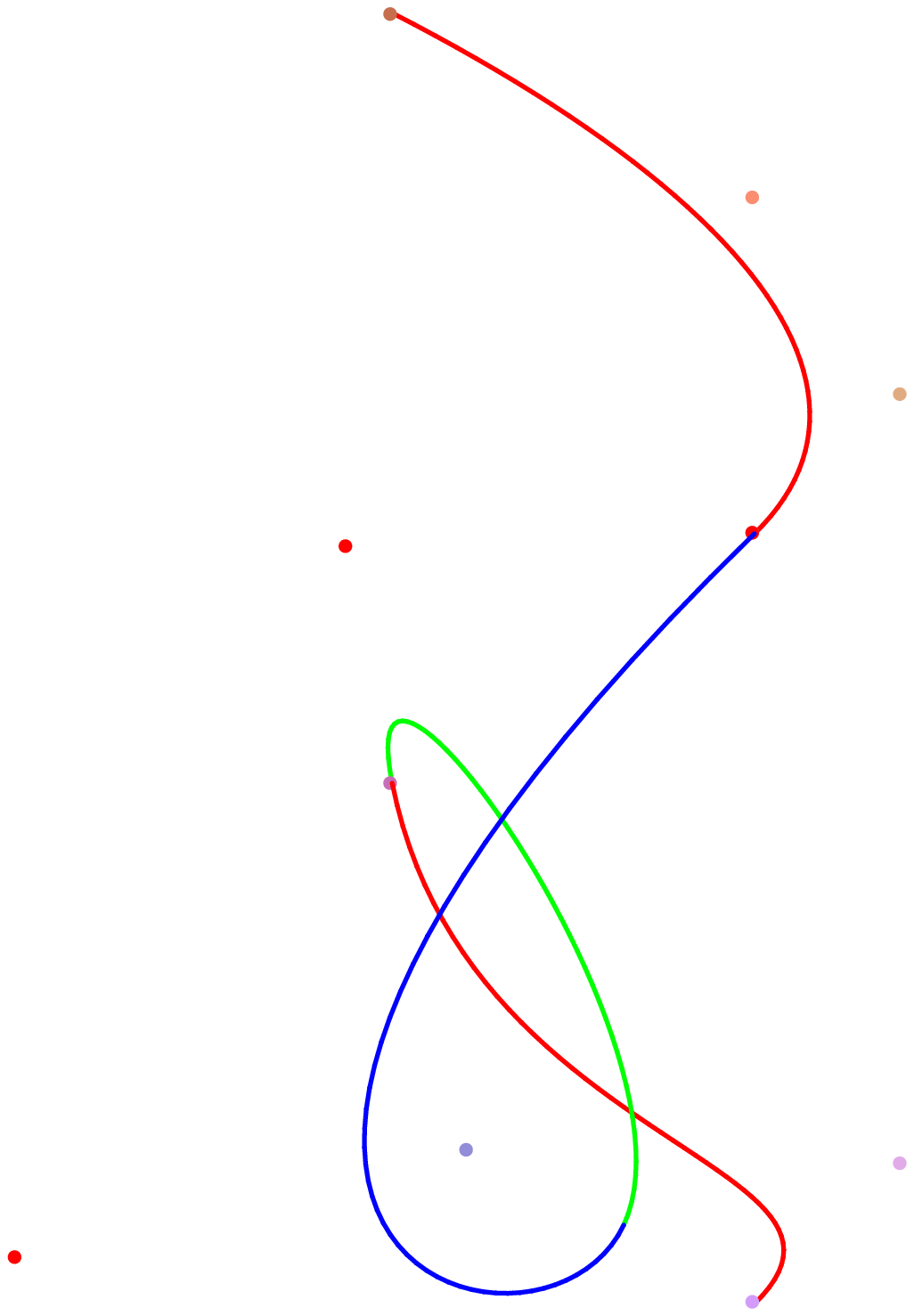}}
\caption{Left shows the solution of the occluded curve as a two-piece cubic space curve in Example \ref{example3}, red curves show the input curves $f_1$ and $f_2$, blue curve shows $f_{01}(t),\;t\in[0,1/2]$ and green curve shows $f_{02}(t),\;t\in[1/2,1]$; right shows the solution of the occluded curve as a two-piece cubic space curve in Example \ref{example4}, red curves show the input curves $f_1$ and $f_2$, blue curve shows $f_{01}(t),\;t\in[0,1/2]$ and green curve shows $f_{02}(t),\;t\in[1/2,1]$.}\label{space1}
\end{figure}

In order to gain more local control and admit more complexity of the occluded curve, sometimes we need to use multi-piece B-spline curve to solve the occluded curve problem. As for the one-piece plane curve, the key is to find an appropriate Lagrangian. Now let us see some examples.
\begin{example}\label{example3}
let a two-piece cubic space curve $f_{01}(t)=\sum_{i=1}^4b_{i,3}^1(t)p_i$ and $f_{02}(t)=\sum_{i=2}^5b_{i,3}^2(t)p_i$ be the solution to the occluded curve which connects two one-piece cubic curves $f_1(t)=\sum_{i=-2}^1b_{i,3}(t)p_i$ and $f_2(t)=\sum_{i=5}^8b_{i,3}(t)p_i$, where $p_{-2}=(0,0,-2),p_{-1}=(1,4,-1.5),p_0=(2,1,-1),p_1=(4,3,-0.5),p_5=(7,2,0),p_6=(8,3,0.5),p_7=(9,1,1),p_8=\\(10,2,1.5)$. The basis functions $b_{i,3}^1(t),t\in[0,1/2]$ and $b_{i,3}^2(t),t\in[1/2,1]$ are defined over the knot sequence $\{0,0,0,0,1/2,1,1,1,1\}$. The Lagrangian we chose in this example is $L=(I_{1,2}\times I_{2,2})\cdot I_{3,2}+I_{1,3}\cdot I_{2,3}$. The solution is shown in Figure \ref{space1}
\end{example}
\begin{example}\label{example4}
Same as in Example \ref{example3}, let a two-piece cubic space curve $f_{01}(t)$ and $f_{02}(t)$ be the solution to the occluded curve. But let the input data now be $p_{-2}=(0,0,0),p_{-1}=(2,1,1),p_0=(1,2,1.5),p_1=(0.5,1.5,2),p_5=(0,0,3),p_6=(-2,1,4),p_7=(1,2,4.5),p_8=\\(0.5,1.5,5)$. The Lagrangian we choose in this example is same as in Example \ref{example3}. The solution is shown in Figure \ref{space1}
\end{example}
Now let us generalize an algorithm of solving the occluded curve problem, which is approximated by uniform B-spline curves that satisfies the Bessel end conditions:
\begin{algorithm}\label{algo1}

\textbf{Input:} Two sets of control points $(p_i)_{i=i_0}^{i_f}$ and $(p_j)_{j=j_0}^{j_f}$, the degree of the occluded curve $d$, the number of the segments of the occluded curve $l$, and a discrete Lagrangian $L$.

\textbf{Output:} A set of control points $(p_k)_{k=i_f}^{j_0}$ which satisfy the extremization of the variational problem subjects to the discrete Lagrangian $L$;

\textbf{Step 1:} Transform the input curves to the normalized position;

\textbf{Step 2:} Put the corresponding differences of control points into $L$, use the method introduced in (\ref{discrete}) to calculate the Euler-Lagrange equations and obtain the solution;

\textbf{Step 3:} Transform the input curves and the solution curves from the normalized position to the original position.
\end{algorithm}
\begin{remark}
Algorithm \ref{algo1} shows us, that for any discrete input data which describe an occluded curve problem approximated by B-spline curve, we can find an aesthetically pleasing solution curve, which extremizes the variational problem subjects to an appropriate Lagrangian.

The choice of the Lagrangian is crucial, not only because it yields different results, but it also decide the solvability of the discrete system. We need to use the Lagrangian to derive sufficient number of Euler-Lagrange equations. For example, using $L=I_{1,2}\cdot I_{2,2}$ will derive two independent Euler-Lagrange equations, while using $L=I_{1,2}+I_{2,2}$ will only derive one Euler-Lagrange equation.

The input data we choose are relatively simple, so we only care about the boundary conditions described in Figure \ref{path1}. But what if we have relatively complex input data, will Algorithm \ref{algo1} still offer us satisfied solutions?
\end{remark}
\section{Adding Complexities to Solution Curve}
\begin{figure}[H]
\begin{pspicture}(-2,3.3)
\psbezier[linewidth=1pt,showpoints=false]{}(2.5,2)(3,4)(3.5,1)(4,3)
\psbezier[linewidth=1pt,showpoints=false]{}(1,1)(1.5,6)(2,0)(2.5,2)
\psbezier[linewidth=1pt,showpoints=false]{}(7,2)(8,3)(9,1)(10,2)
\pscircle[linewidth=1pt,linestyle=dotted](5.5,2.5){1.58}
\rput(1.8,2){$f_1$}
\rput(8,2){$f_2$}
\rput(5.5,2.5){$?$}
\end{pspicture}
  \caption{Adding complexity to a path connecting the paths $f_1$ and $f_2$}\label{path_mul}
\end{figure}
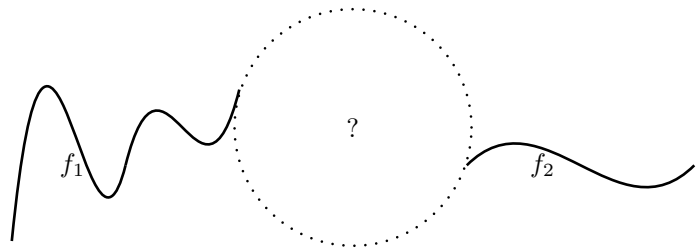

From the construction of the occluded curve we can see, that the complexity of the solution curve highly depends on the input data, the degree of the solution curve, and the order of difference invariants in Lagrangian. Using one piece cubic spline as solution curve limits its complexity. As shown in the previous examples, one piece cubic spline curve will be enough to solve the occluded curve problem described in Figure \ref{path}, but what should the solution look like in Figure \ref{path_mul}?
\subsection{Complexity of Occluded Curve}
The solution curve to the problem described in Figure \ref{path_mul} requires complexity with more than one inflection point. It can be obtained by inserting control points, which will result in either increasing the degree of the solution curve and keep the number of spline piece; or increasing the number of spline piece and keep the degree of the solution curve. Let us now discuss the case of increasing complexity of the solution curve by inserting one control point.

\begin{figure}[H]
  \centerline{\includegraphics[height=4cm]{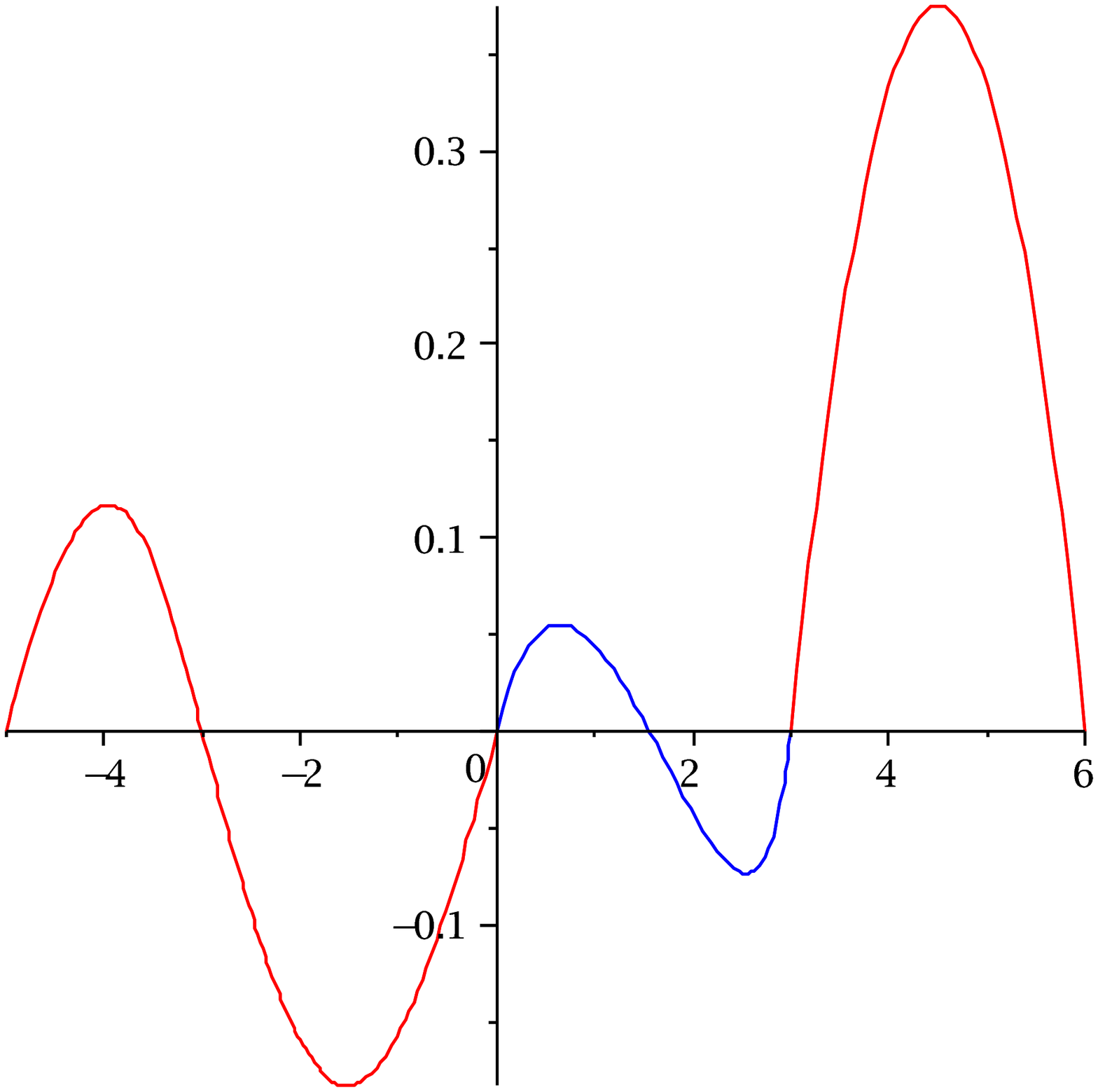}\qquad\includegraphics[height=4cm]{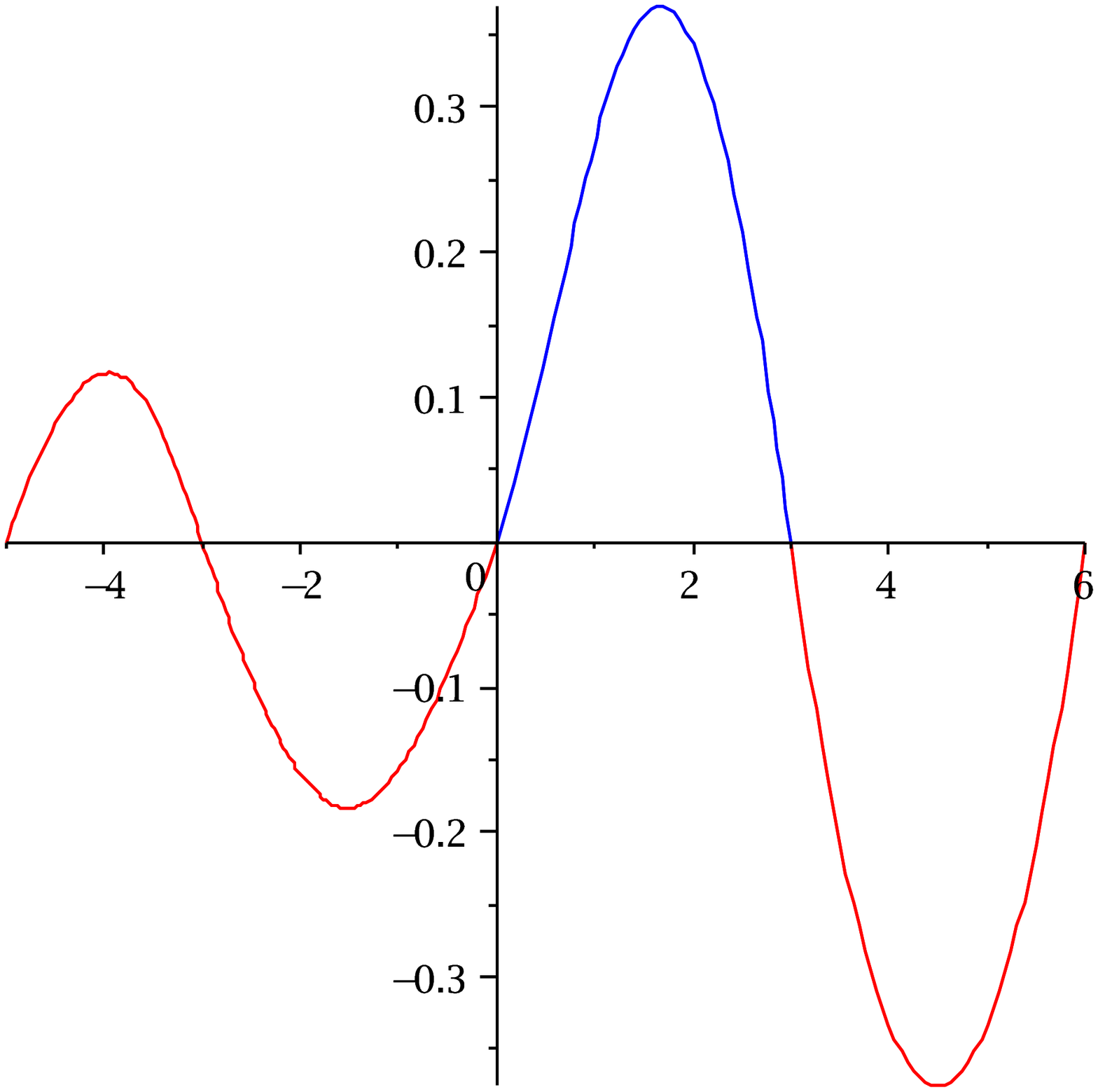}}
\caption{Left shows a one piece cubic solution curve with one inflection point, corresponding to the same signs of the tangent vectors at the entering and leaving points of the occluded curve; right shows a one piece cubic solution curve with no inflection point, corresponding to the different signs of the tangent vectors. The Lagrangians used in both cases are the same as in Example \ref{example1}.}\label{mul_00}
\end{figure}

The complexity of a one piece cubic solution curve will be determined by the input data, and can be classified into two cases:
\begin{itemize}
\item Case 1: The tangent vectors at the entering (left) and leaving (right) point of the occluded curve have same signs. In this case, the solution curve has one inflection point. An example is shown in the left picture of Figure \ref{mul_00}.
\item Case 2: The tangent vectors at the entering (left) and leaving (right) point of the occluded curve have different signs. In this case, the solution curve has no inflection point. An example is shown in the right picture of Figure \ref{mul_00}.
\end{itemize}
\subsection{Adding Complexity}
Adding one more control point to a one piece cubic solution curve will either make it become a two piece cubic curve or an one piece quartic curve. As we can see in the Figure \ref{mul_00}, adding one control point to the occluded curve in the left picture will not add inflection point; adding one control point to the right picture will add two more inflection points to the occluded curve. We need to discuss these two cases respectively.
\subsubsection{Case 1}
Let an one piece cubic solution curve be $f(t)=\sum_{i=1}^4b_i(t)p_i$, where the basis functions are defined over $(0,0,0,0,1,1,1,1)$, after inserting one control point, it becomes either of the following curves.

\textbf{Two piece cubic curve:} $f_j(t)=\sum_{i=j}^{j+3}{b}^j_{i-j+1}(t)q_i,j=1,2$, where the basis functions are defined over $(0,0,0,0,1/2,1,1,1,1)$;

\textbf{One piece quartic curve:} $f(t)=\sum_{i=1}^5b_i(t)q_i$, where the basis functions are defined over $(0,0,0,0,0,1,1,1,1,1)$.

Using either two piece cubic curve or one piece quartic curve as solution curve we will have five control points in total $q_i=(x_i,y_i),\,i=1..5$, which can construct three second order difference invariants. Hence we can use the Lagrangian: $L=(I_{1,2}\cdot I_{2,2})(I_{2,2}\cdot I_{3,2})$, where $I_{i,k}=q_{i+1}^{k-1}-q_{i}^{k-1})$ are difference invariants at the normalized position, where the last control point of the entering (left) curve at the origin and the first control point of leaving (right) curve on the $x$-axis.

From the boundary conditions in (\ref{boundary}), we obtain two unknowns $\alpha$ and $\beta$; other two unknowns are the coordinates of $q_3$. Since we only have three Euler-Lagrange equations, we need to use a boundary condition to eliminate one unknown in $q_3$. Judging from the normalized position, noticing that in this case inserting control point $q_3$ will not add inflection point, so we do not need to constrain the $y$-coordinate of $q_3$. Therefore, set $x_3=(x_2+x_4)/2$,  solve the three Euler-Lagrange equations for $\alpha$, $\beta$ and $y_3$, we will obtain the solution curve. Next let us see an example.
\begin{example}\label{mul_0_1}
Input are two cubic curves, one is defined over $(0,0,0,0,1/4,1/2,3/4,1,1,1,1)$ with control points $(-6,0),(-5,-1),(-4,0),(-3,1),(-2,0),(-1,-1),(0,0)$; and another is defined over $(0,0,0,0,1,1,1,1)$ with control points $(3,0),(4,2),(5,2),(6,0)$. A two piece cubic solution curve is defined over $(0,0,0,0,1/2,1,1,1,1)$ and an one piece quartic solution curve is defined over $(0,0,0,0,0,1,1,1,1,1)$. Using the Lagrangian $L=(I_{1,2}\cdot I_{2,2})(I_{2,2}\cdot I_{3,2})$, the solution curves are shown in Figure \ref{mul_01}
\begin{figure}[H]
  \centerline{\includegraphics[height=4cm]{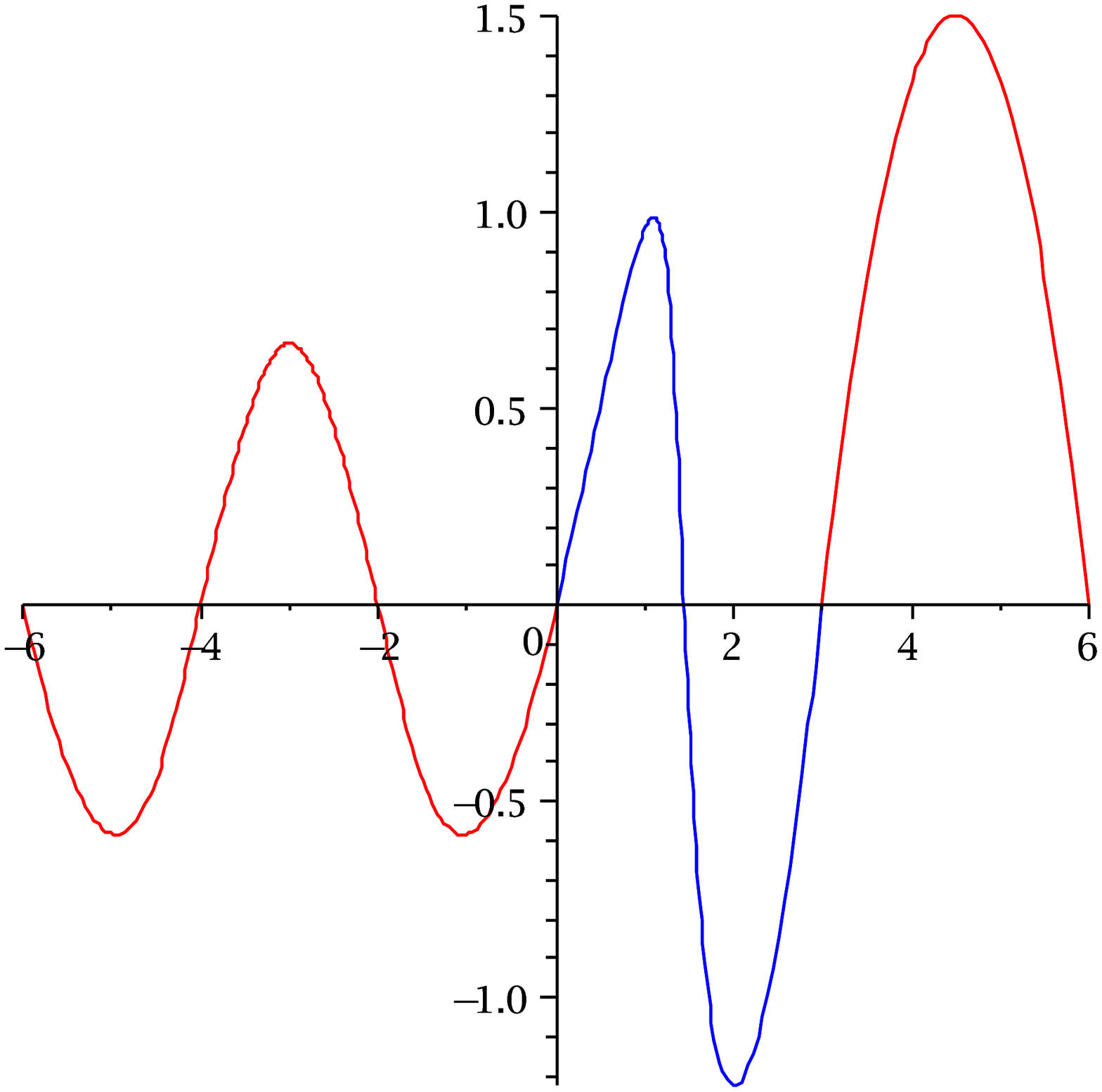}\qquad\includegraphics[height=4cm]{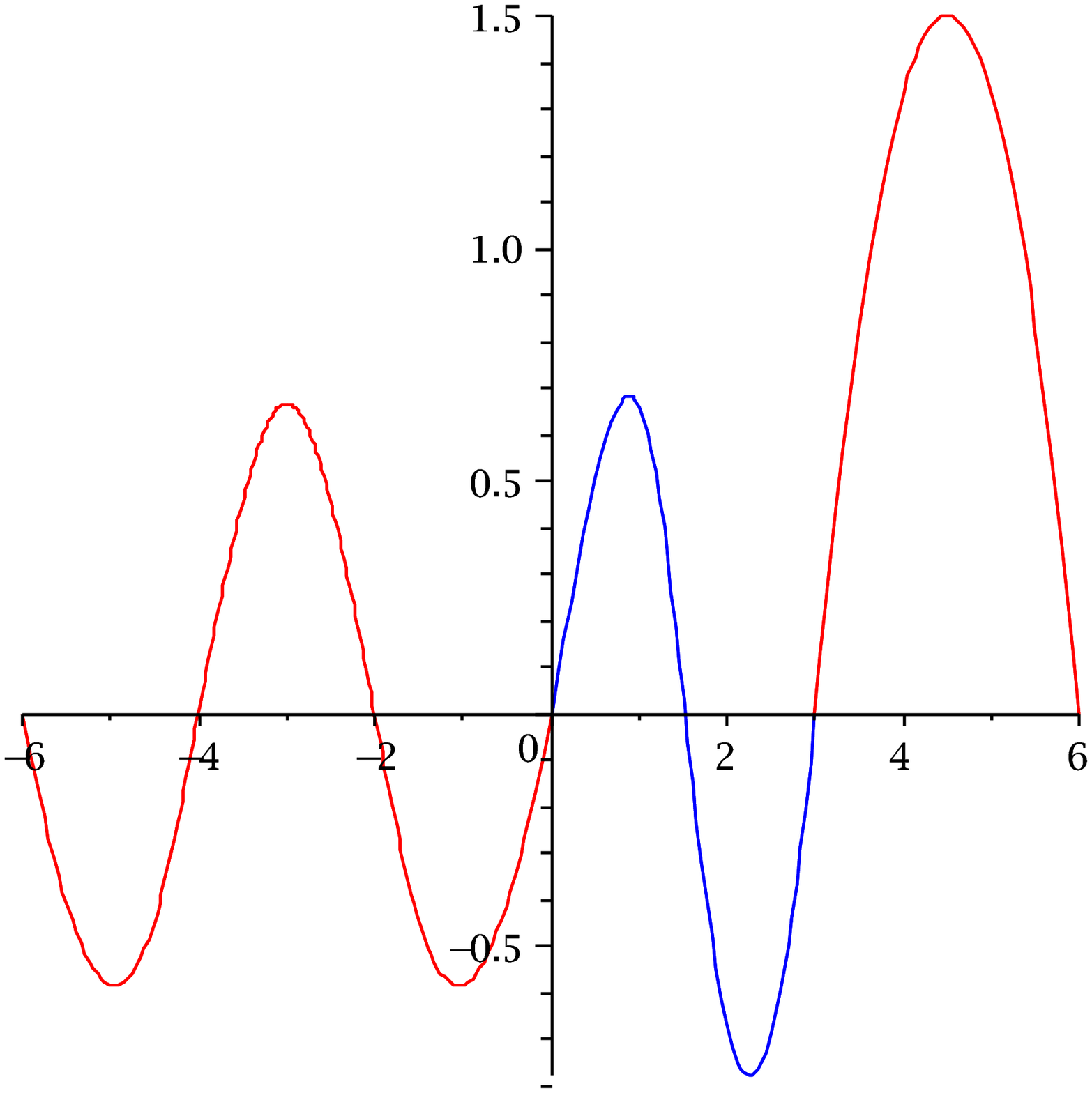}}
\caption{The solution curves of Example \ref{mul_0_1}. Left shows the solution curve as a two piece cubic spline curve, whose basis functions are defined over $(0,0,0,0,1/2,1,1,1,1)$ at the normalized positions; right shows the solution curve as an one piece cubic spline curve, whose basis functions are defined over $(0,0,0,0,0,1,1,1,1,1)$ at normalized positions.}\label{mul_01}
\end{figure}
\end{example}
\subsubsection{Case 2}
All setups are the same as in case 1, except the constrain to eliminate one unknown in $q_3$. Since we want to add inflection points to the solution curve, which means $y_3$ need to have different sign as $y_2$ and $y_4$. Set $y_3=-(y_2+y_4)/2$, solve the three Euler-Lagrange equations for $\alpha,\beta$, and $x_3$, we obtain the solution curve. Let us see an example:
\begin{example}\label{mul_1_1}
Input are two cubic curves, one is defined over $(0,0,0,0,1/3,2/3,1,1,1,1)$ with control points $(-2,6),(-1,3),(0,6),(1,4),(2,3),(4,3)$; and another is defined over $(0,0,0,0,1,1,1,1)$ with control points $(7,2),(8,1),(9,3),(10,2)$. A two piece cubic solution curve is defined over $(0,0,0,0,1/2,1,1,1,1)$ and an one piece quartic solution curve is defined over $(0,0,0,0,0,1,1,1,1,1)$. Using the same Lagrangian as in Example \ref{mul_0_1}, the solution curves are shown in Figure \ref{mul_11}
\begin{figure}[H]
  \centerline{\includegraphics[width=4cm]{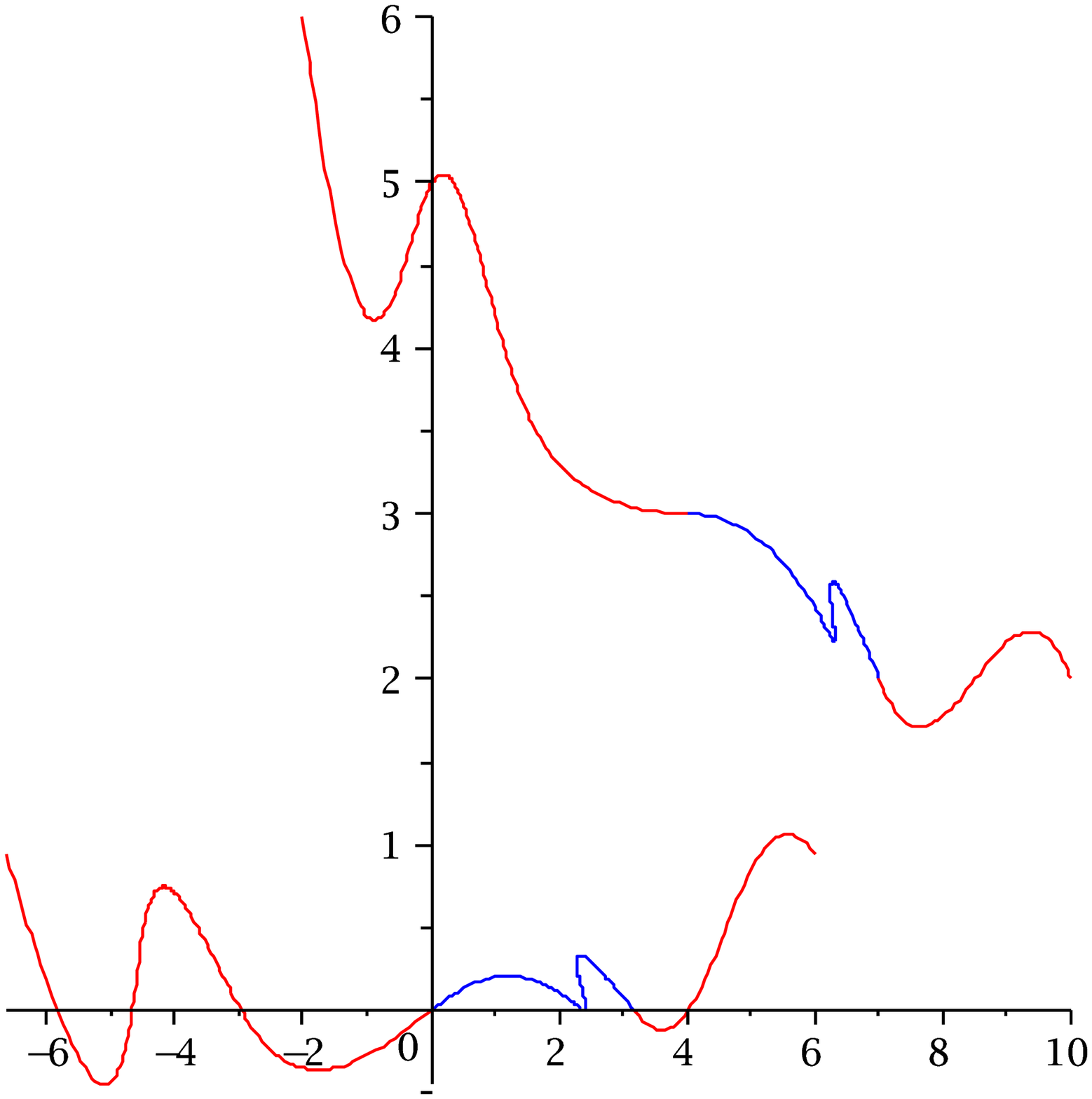}\includegraphics[width=4cm]{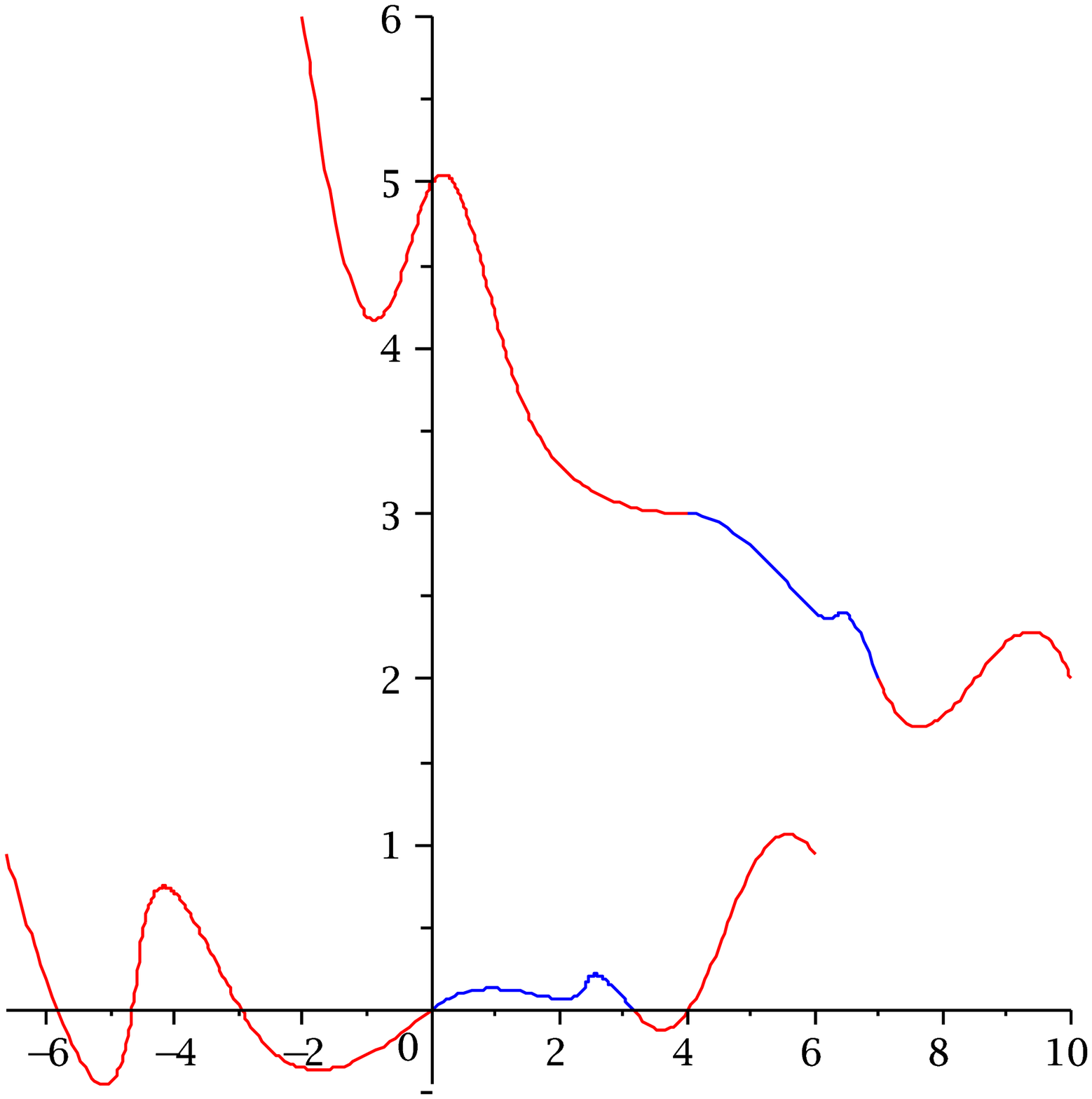}}
\caption{The solution curves of Example \ref{mul_1_1}. Left shows the solution curve as a two piece cubic spline curve, whose basis functions are defined over $(0,0,0,0,1/2,1,1,1,1)$ at the original and normalized positions; right shows the solution curve as an one piece quartic spline curve, whose basis functions are defined over $(0,0,0,0,0,1,1,1,1,1)$ at the original and normalized positions.}\label{mul_11}
\end{figure}
\end{example}
\subsection{Adding More Constrains}
\begin{figure}[H]
  \centerline{\includegraphics[width=5cm]{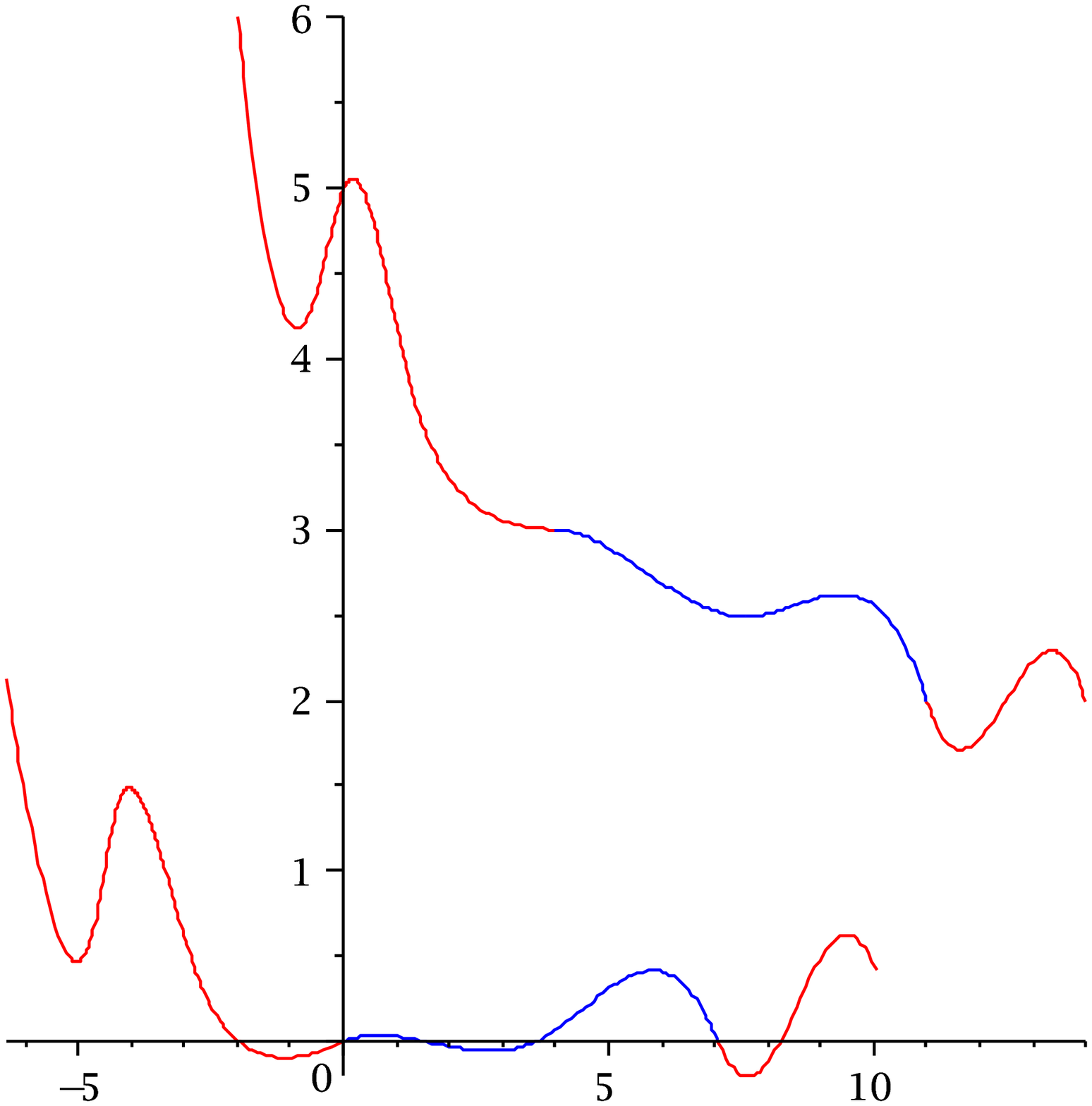}\qquad\includegraphics[width=5cm]{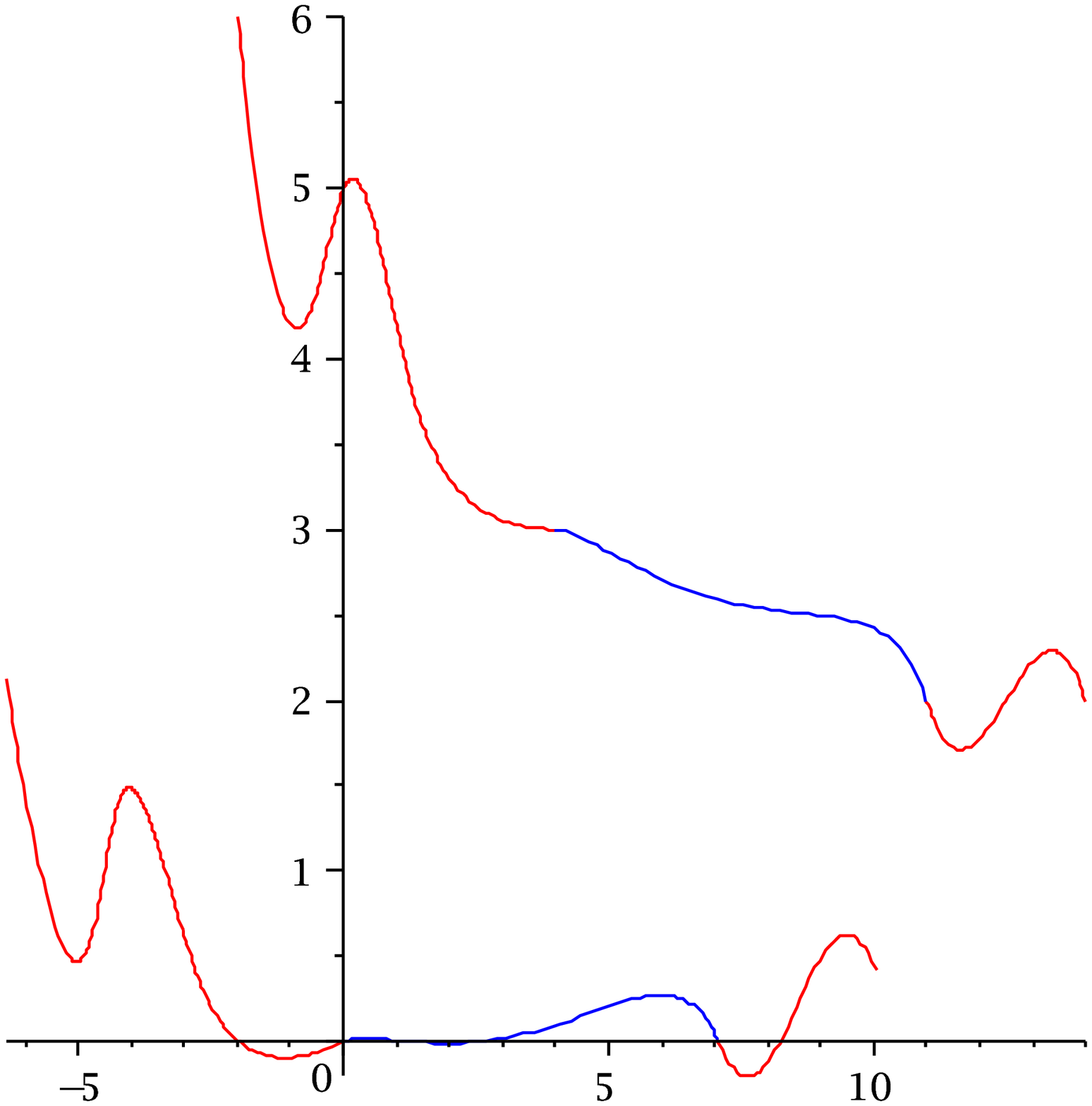}}
\caption{Solution to Example \ref{mul_1_1}, comparing with Figure \ref{mul_11}, the gap between two input curves are increased. }\label{mul_21}
\end{figure}

From the example in case 1 we can see that quartic solution curve looks more smooth than the cubic solution curve; from the solution in case 2 we can see, after two inflection points are added, the solution curves seem a little too ``complex'', the curvature around the inflection point is too high. Hence we need to consider the length of the interval between the starting and ending points of the solution curve on $x$-axis. If we increase gap, the solution curve will become more smooth, as shown in Figure \ref{mul_21}.

\begin{figure}[H]
  \centerline{\includegraphics[width=5cm]{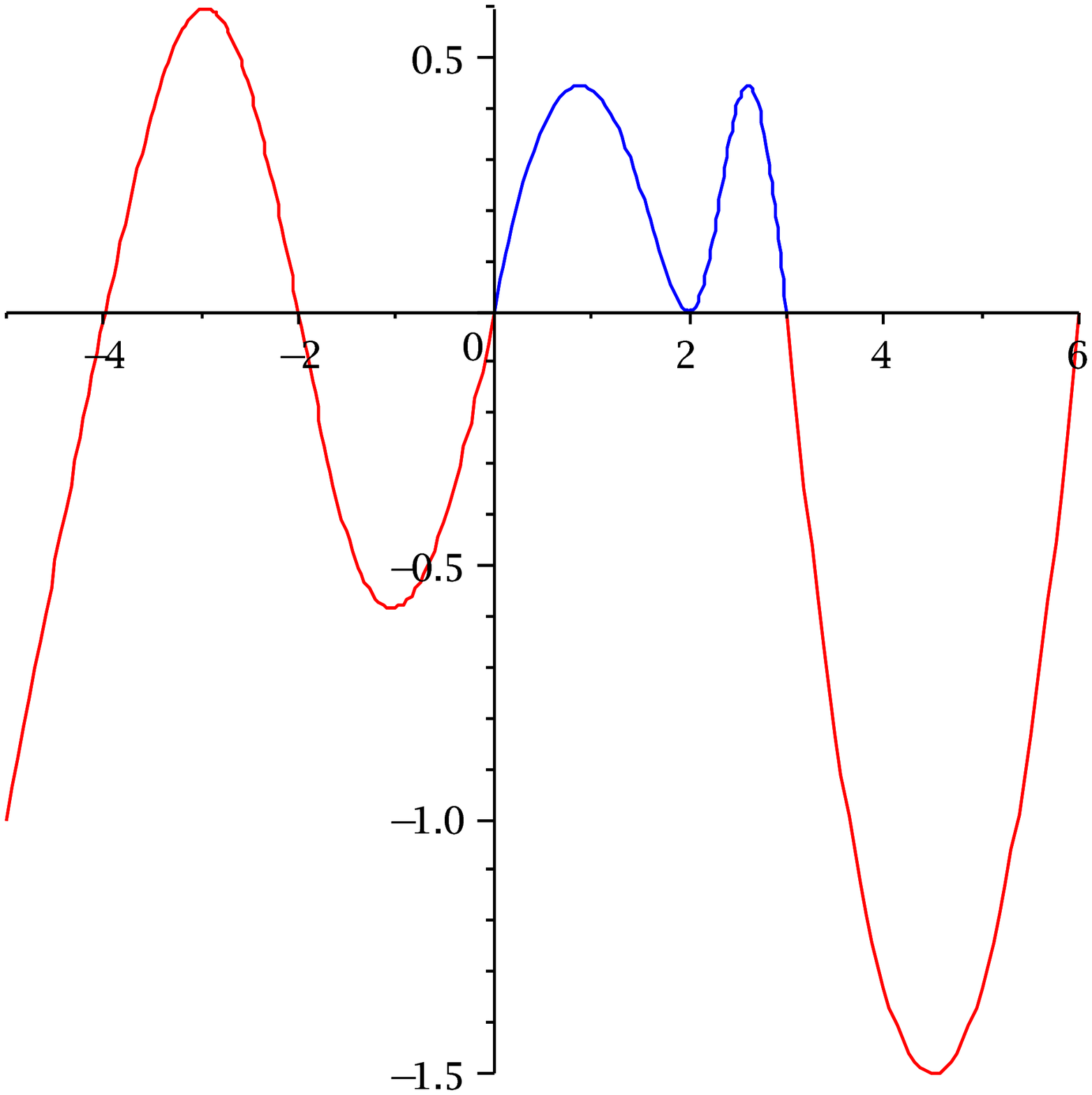}\qquad\includegraphics[width=5cm]{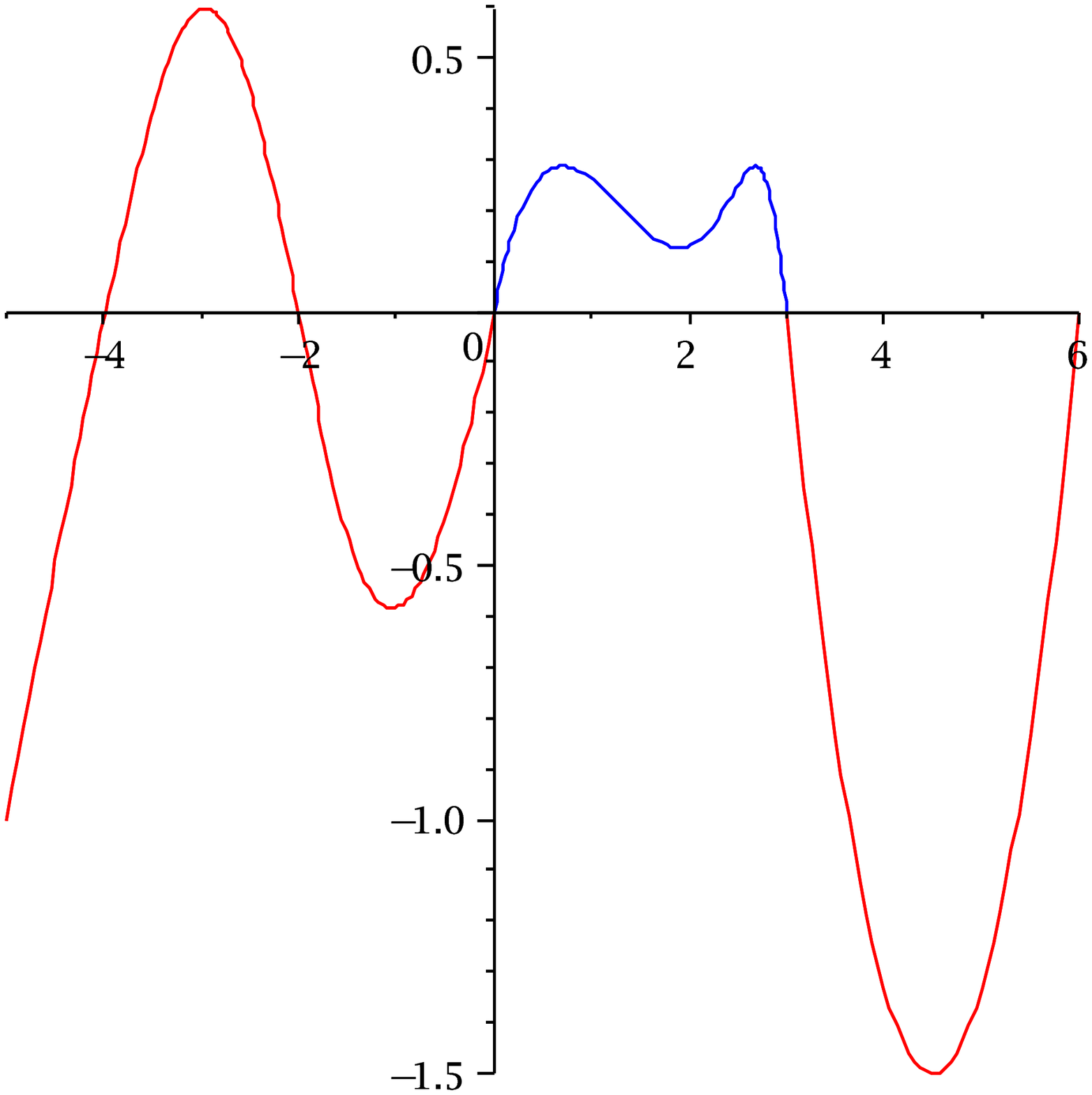}}
\caption{Left shows a two piece cubic solution curve, right shows an one piece quartic solution curve}\label{mul_31}
\end{figure}

\begin{figure}[H]
  \centerline{\includegraphics[width=5cm]{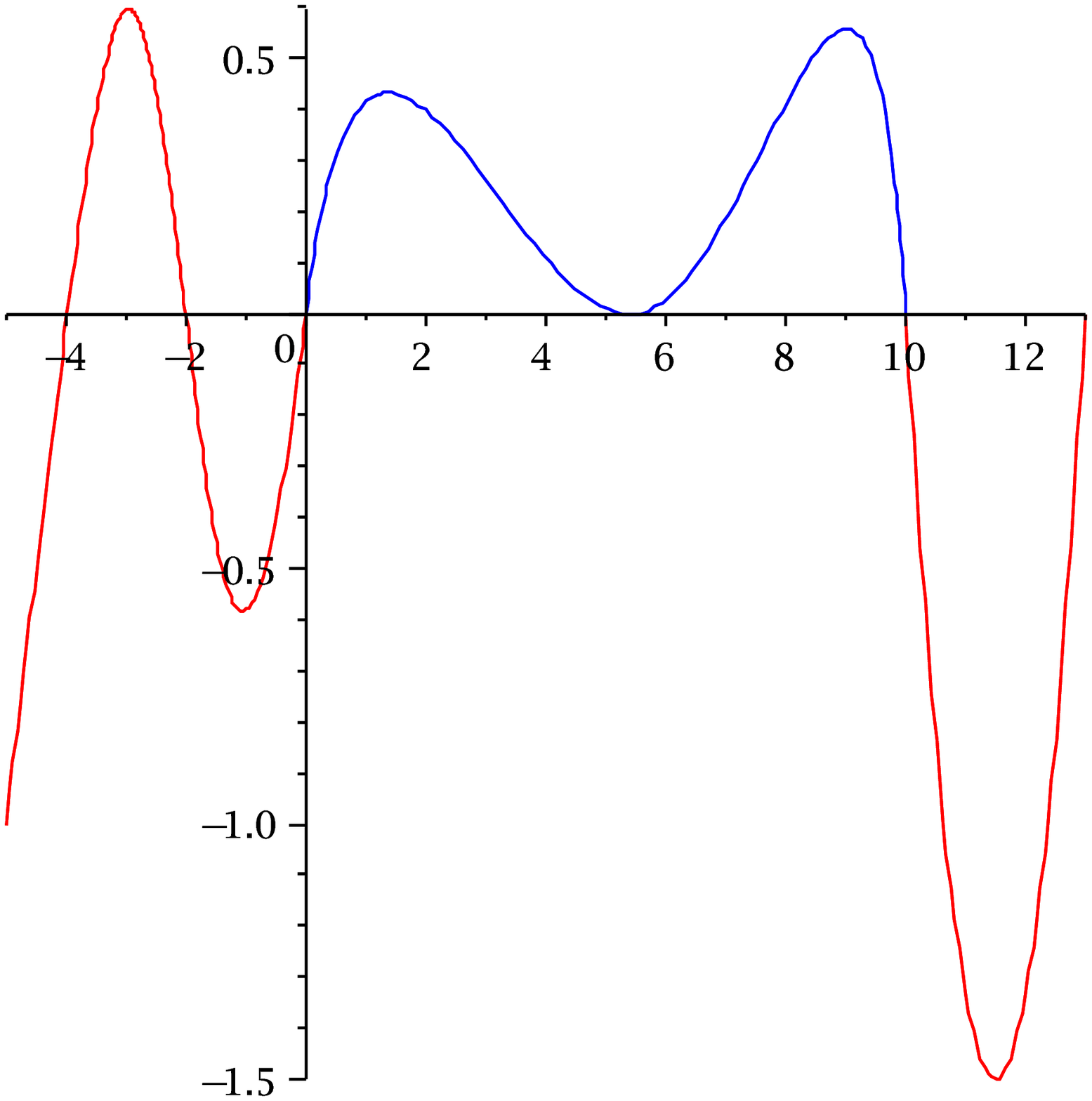}\qquad\includegraphics[width=5cm]{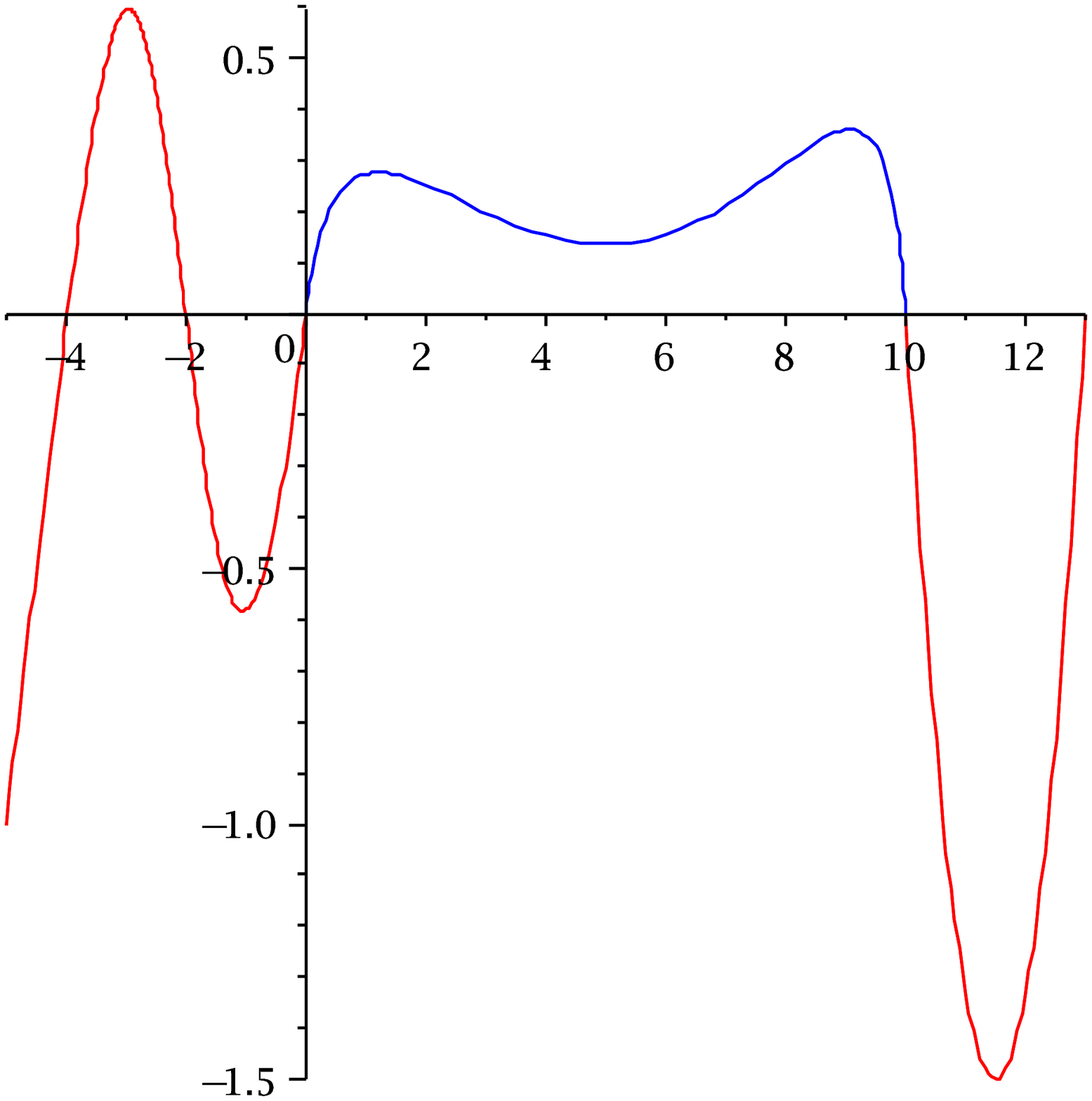}}
\caption{Using same input data as in Figure \ref{mul_31}, but increase the gap between the two given curves. Left shows a two piece cubic solution curve, right shows an one piece quartic solution curve}\label{mul_41}
\end{figure}

Solution curves of other examples which using the same Lagrangian as in Example \ref{mul_01}, using one piece quartic and two piece cubic spline curves, with different gaps are shown in Figure \ref{mul_31} and \ref{mul_41}.
Assuming that the left input curve has more complexity than the right input curve, now we can generalize a revised algorithm:
\begin{algorithm}

\textbf{Input:} Two sets of control points $(p_i)_{i=i_0}^{i_f}$ and $(p_j)_{j=j_0}^{j_f}$, the degree of the occluded curve $d$, the number of the segments of the occluded curve $l$, and a discrete Lagrangian $L$.

\textbf{Output:} A set of control points $(p_k)_{k=i_f}^{j_0}$ which satisfy the extremization of the variational problem subject to the discrete Lagrangian $L$;

\textbf{Step 1:} Transform the input curves to the normalized position;

\textbf{Step 2:} Use the gap between two input curves as a measure, calculate the number of the inflection points $n_1$ and $n_2$ in the same gap length of the left and right input curves which are adjacent to the occluded curve. The number of the inflection points for the solution curve is $n=(n_1+n_2)/2$, and rounding to the integer part. Depending on the degree of the solution curve $d$, determine the number of spline pieces. Using a one piece cubic curve as basic solution curve, adding control points to it.

\textbf{Step 3:} Judging from different cases introduced previously, using boundary conditions and constrains to eliminate some unknowns, put the difference invariants into $L$. Use the method introduced in (\ref{discrete}) to calculate the Euler-Lagrange equations and obtain the solution;

\textbf{Step 4:} Transform the input curves and the solution curves from the normalized position to the original position.
\end{algorithm}
\section{Conclusion}
We have showed that a complex variational problem in the smooth case can be simplified by using an approximation method, such as B-spline curves, and B-spline approximated occluded curve can be obtained with relative ease and aesthetically pleasing effects.

When complexity is not required, the solution curve can be approximated by a one piece cubic B-spline curve. Higher order B-spline curves are used when more smoothness is required.

When using a multi-piece B-spline curve as solution curve, we need to find an appropriate Lagrangian which can derive a sufficient number of Euler-Lagrange equations to calculate the unknowns. For example, a Lagrangian such as $L=I_{1,2}+I_{2,2}$ will not give us two independent Euler-Lagrange equations. When the number of unknowns exceeds the number of Euler-Lagrange equations we can obtain, we need to add more boundary conditions.

When more complexity is required, a more detailed solution curve can be constructed by adding more control points to the original one piece cubic curve. The result will be a multi-piece cubic curve or a one piece higher order curve. The additional unknowns can be solved for using extra boundary conditions in order to obtain more inflection points in the solution curve.

Transforming the input curves from the original position to a normalized position is not necessary, but will help us to simplify the problem. At the normalized position, the complexity of the solution curve can be easily analyzed. The normalization process is a simple example of a general procedure called ``moving frames'', which can be used for general Lie group actions.\cite{mansfield}\cite{olver}\cite{felsnolver}\cite{hubert}.

Some extreme cases are not considered in this paper, since the solution are trivial. For example, when both input curves have zero slope at the boundary points, the solution is obviously a straight line.

Compared with using a standard numerical solution to a smooth Euclidean invariant variational problem\cite{mansfield}, the method used here offers more computational advantage, flexibility and reliability.

Further study could focus on the occluded surface problem and discrete variational problems other than occluded curve problem. Casting the variational problem in Conformal Geometric Algebra will be interesting, challenging, and could lead to further progress.

\textsl{This paper is supported by EPSRC grant EP/E001823/1.} 
\bibliographystyle{plain}
\bibliography{paper}
\end{document}